\documentclass[11pt]{amsart}
\usepackage{amsfonts}
\usepackage{amssymb}
\usepackage[letterpaper, margin=0.6in]{geometry}
\usepackage{fancyhdr}

\input{tcilatex}

\begin{document}
\date{}
\author{Ali Enayat}
\author{Albert Visser}
\title[]{Incompleteness of boundedly axiomatizable theories \\
}

\maketitle

\begin{abstract}
\noindent Our main result (Theorem A) shows the incompleteness of any
consistent sequential theory $T$ formulated in a finite language such that $%
T $ is axiomatized by a collection of sentences of bounded
quantifier-alternation-depth. Our proof employs an appropriate reduction
mechanism to rule out the possibility of completeness by simply invoking
Tarski's Undefinability of Truth theorem. We also use the proof strategy of
Theorem A to obtain other incompleteness results (as in Theorems A$^{+},$ B
and B$^{+}).$
\end{abstract}

\bigskip

\noindent Our main result (Theorem A) was prompted by the following question%
\footnote{%
Thanks to Roman Kossak for bringing this interesting question to our
attention. We are also indebted to Mateusz \L e\l yk, Dino Rossegger, and
Saeed Salehi for their helpful feedback on the preliminary drafts of this
paper, and to Emil Je\v{r}\'{a}bek for his perceptive remarks on the
penultimate draft.} of Steffen Lempp and Dino Rossegger; the question arose
in the context of their joint work \cite{Uri-Steffen-Dino} with Uri Andrews,
David Gonzalez, and Hongyu Zhu, in which they establish: \textit{For a
complete first-order theory }$T$\textit{, the set of models of }$T$\textit{\
is }$\Pi _{\omega }^{0}$\textit{-complete under Wadge reducibility (i.e.,
reducibility via continuous functions) if and only if }$T$\textit{\ does not
admit a first-order axiomatization by formulae of bounded quantifier
complexity}. In what follows $\mathsf{PA}^{-}$ is the well-known finitely
axiomatized fragment of $\mathsf{PA}$ (Peano Arithmetic) whose axioms
describe the non-negative substructure of discretely ordered rings (with no
instance of the induction scheme, hence the minus superscript), as in Kaye's
text \cite{Kaye's text} on models of $\mathsf{PA}$. \medskip

\noindent \textbf{Question }$\diamondsuit $.~\textit{Is there a consistent
completion of }$\mathsf{PA}^{-}$ \textit{that is axiomatized by a set of
sentences of bounded quantifier complexity}?\medskip

\noindent \textbf{Remark 1.}~It is well-known that the answer to the above
question is in the negative when $\mathsf{PA}^{-}$ is strengthened to $%
\mathsf{PA}$. This result follows from a theorem of Rabin \cite{Rabin
Arithmetic} that states that for each $n\in \omega $ no consistent extension
of $\mathsf{PA}$ (in the same language) is axiomatized by a set of $\Sigma
_{n}$-sentences.\footnote{%
In his paper Rabin points out that this result was possibly known (but not
published) by others, including Feferman, Wang, Scott, Kreisel, and
Tennenbaum.} Rabin's result refines an earlier theorem of Ryll-Nardzewski 
\cite{Ryll-Nardzewski} that states that no consistent extension of $\mathsf{%
PA}$ is finitely axiomatizable. Ryll-Nardzewski and Rabin both employed
model-theoretic arguments relying on nonstandard elements to prove the
aforementioned results (see Theorem 10.2 on p.132 of Kaye's text \cite%
{Kaye's text} for a modern treatment). Rabin's result can be also
established with an argument that mixes proof-theoretic machinery, partial
satisfaction classes, and G\"{o}del's second incompleteness theorem (see
Theorem 2.36 of Chapter III\ of \cite{Hajek-Pudlak} for an exposition). As
shown by Montague \cite{Montague-sem-closed}, a similar result can be
established for any \textit{inductive} sequential theory $T$, i.e., a
sequential theory that has the power to prove the full scheme of induction
over its `natural numbers' for all formulae in the language of $T$ (this
involves a designated interpretation of a suitable base theory of the
natural numbers). In this more general setting the relevant hierarchy is
based on the depth of quantifier alternations; canonical examples of
inductive sequential theories include all extensions of $\mathsf{PA}$, $%
\mathsf{Z}$ (Zermelo set theory), $\mathsf{Z}_{2}$ (second order
arithmetic), and $\mathsf{KM}$ (Kelley-Morse theory of classes).\medskip

Theorems A and A$^{+}$ are formulated for sequential theories. At first
approximation, a theory is sequential if it supports a modicum of coding
machinery to handle finite sequences of all objects in the domain of
discourse. Sequentiality is a modest demand for theories of arithmetic and
set theory; however, by a theorem of Visser \cite{Visser on Q}, (Robinson's) 
$\mathsf{Q}$ is not sequential. There are many equivalent definitions of
sequentiality; the original definition due to Pudl\'{a}k (used by Je\v{r}%
\'{a}bek \cite{JerabekPAminus} in his proof of sequentiality of $\mathsf{PA}%
^{-}$) is as follows: A theory $T$ is sequential if there is a formula $N(x)$%
, together with appropriate formulae providing interpretations of equality,
and the operations of successor, addition, and multiplication for elements
satisfying $N(x)$ such that $T$ proves the translations of the axioms of $%
\mathsf{Q}$ (Robinson's arithmetic) when relativized to $N(x)$; and
additionally, there is a formula $\beta (x,i,w)$ (whose intended meaning is
that $x$ is the $i$-th element of a sequence $w$) such that $T$ proves that
every sequence can be extended by any given element of the domain of
discourse, i.e., $T$ proves:

\begin{center}
$\forall w,x,k\ \exists w^{\prime }\ \forall i,y\ \left[ \left[ N(k)\wedge
i\leq k\right] \rightarrow \left[ 
\begin{array}{c}
\beta (y,i,w^{\prime })\leftrightarrow \\ 
\left[ i<k\wedge \beta (y,i,w)\right] \vee \left[ i=k\wedge y=x\right]%
\end{array}%
\right] \right] .$
\end{center}

\noindent For more information about sequentiality, see \cite{Visser-Small}%
.\medskip

\noindent In light of the aforementioned proof\footnote{%
Indeed Je\v{r}\'{a}bek's result is stronger since it establishes the
sequentiality of a weaker theory than the usual formulation of $\mathsf{PA}%
^{-}$. This weaker theory is a universal theory and described by Je\v{r}\'{a}%
bek as the theory of discretely ordered commutative semirings with a least
element, a theory in which the existence of predecessors is unprovable.} of
sequentiality $\mathsf{PA}^{-}$ by Je\v{r}\'{a}bek \cite{JerabekPAminus},
the following general result answers Question $\diamondsuit $ in the
negative. Note the condition of finiteness of the language in Theorem A
cannot be eliminated, as indicated in Remark 4.

\begin{itemize}
\item Throughout the paper $\Sigma _{n}^{\ast }$ is the hierarchy of
formulae whose measure of complexity is \textit{depth of quantifier
alternation}, as in \cite{Visser-Small}.\medskip
\end{itemize}

\noindent \textbf{Theorem A.}~\textit{For any fixed} $n\in \omega ,$ \textit{%
every} \textit{consistent sequential theory formulated in a finite language
that is axiomatized by a set of }$\Sigma _{n}^{\ast }$-\textit{sentences is
incomplete}.\medskip

\noindent \textbf{Remark 2.}~The proofs of our other theorems are all based
on the proof of Theorem A; we provide a direct proof of this theorem,
employing a middle range level of abstraction, addressed to a non-specialist
reader. The three ingredients of the proof are (1) partial satisfaction
predicates, (2) Rosser's `trick' used in his celebrated generalization of G%
\"{o}del's first incompleteness theorem in which the technical hypothesis of 
$\omega $-consistency (more precisely: 1-consistency) is removed; and (3)
Tarski's Undefinability of Truth theorem\footnote{%
Tarski's theorem \cite[p.46]{Tarski-Mostwoski-Robinson} is very general; it
states that if $T$\ is a theory formulated in a language $\mathcal{L}$ has
the property that the diagonal function $\varphi (x)\longmapsto \varphi
(\ulcorner \varphi \urcorner )$ is representable in $T$, then there is no $%
\mathcal{L}$-formula $V(x)$ such that $T$ proves $\psi \leftrightarrow
V(\ulcorner \psi \urcorner )$ for all $\mathcal{L}$-sentences $\psi $. It is
well-known that the diagonal function is representable in the theory $%
\mathsf{R}$ of \cite{Tarski-Mostwoski-Robinson}; which in turn makes it
clear that the diagonal function is representable in sequential theories
formulated in a finite language since (Robinson's) $\mathsf{Q}$, and a
fortiori $\mathsf{R}$, is interpretable in sequential theories.}. Partial
satisfaction predicates were first introduced by Mostowski \cite%
{Mostowski-PAisREF} in his proof demonstrating that $\mathsf{PA}$ is
essentially reflexive (i.e., $\mathsf{PA}$ and all of its extensions in the
same language can prove the formal consistency of each of their finitely
axiomatized subtheories). The technology of partial satisfaction predicates
was further elaborated in the context of sequential theories by Montague 
\cite{Montague-sem-closed}, Pudl\'{a}k \cite{Pudlak-survey}, and finally
Visser \cite{Visser-Small}, whose variant we rely upon in our proof. As we
will see, our proof shows that for complete sequential theories formulated
in a finite language that are axiomatized by a set of $\Sigma _{n}^{\ast }$
statements, the following equation holds:

\begin{center}
Rosser provability from True $\Sigma _{n}^{\ast }$ statements = Truth.
\end{center}

\noindent \textit{Thus the only use of diagonalization in our proof is the
rudimentary Tarskian one}. It should be stressed that upon the completion of
our paper we were informed by Emil Je\v{r}\'{a}bek that he struck upon the
same idea in 2016 and established an abstract version of the above equation
in a MathOverflow answer (as explained in Remark 6). In the subsequent
discussions with Emil, he kindly expressed a preference for us to proceed
with the publication of the paper with the minimal proviso that his
aforementioned prior work be referenced, which is precisely the course of
action the authors settled on.\medskip

There is also a conceptual/pedagogical take-away to our approach. Let $%
\mathsf{R}$ be the well-known fragment of $\mathsf{PA}$ introduced in the
Tarski-Mostowski-Robinson monograph \cite{Tarski-Mostwoski-Robinson} within
which all recursive functions are representable. One can prove the (first)
incompleteness theorem for a consistent computably enumerable extension $T$
of $\mathsf{R}$, without any extra soundness assumptions about $T$, by first
proving Tarski's undefinability of truth theorem with a straightforward
diagonalization (with no need for the fixed point theorem, as in \cite[p.46]%
{Tarski-Mostwoski-Robinson}), and then the incompleteness of $T$ can be
demonstrated using a reductio ad absurdum by verifying that the completeness
of $T$ implies that Rosser provability from $T$ yields a truth definition%
\footnote{%
Using additional machinery, G\"{o}del's second incompleteness theorem can
also be derived from Tarski's Undefinability of Truth theorem; see \cite%
{Visser (G2 from T)}.} (technically, this falls under our Theorem A$^{+}$,
by setting $A=\varnothing $ in that theorem). Note that in contrast to the
usual proof of the incompleteness theorem using the fixed point theorem, our
proof is not constructive, i.e., it does not yield an algorithm that takes a
description of a consistent computably enumerable extension $T$ of $\mathsf{R%
}$ as input and outputs a sentence that is independent of $T$. \medskip

\noindent \textbf{Remark 3.}~As an alternative to the proof of Theorem A
presented below, one can also derive Theorem A (using Fact F below) from a
version of Rosser's Theorem due to Saeed Salehi; see \cite{Salehi-reunion}.
The version proposed by Salehi holds under certain abstract conditions. More
specifically, viewed as an application of Salehi's result, our proof
presented gives a realization of these conditions in the case at hand and a
verification that our realization works. Our proof follows an alternative
route to the Rosser-style result since we employ a reduction to Tarski's
undefinability of truth, where Salehi presents an argument that follows the
traditional Rosser argument more directly. Our different route can also be
employed at the level of abstraction of Salehi's work. Also, as kindly
pointed out by Mateusz \L e\l yk, Theorem A follows from Proposition 31 of
his joint recent work \cite{LW-universal} with Bartosz Wcis\l o (asserting
the existence of so-called $(n,k)$-flexible formulae for computably
enumerable sequential theories).\medskip

\noindent Before presenting the proof of Theorem A, we state an important
fact that plays a crucial role in the proof of Theorem A. The following
result was established by Visser in \cite{Unprov-small-inconsistency}, and
refined in \cite{Visser-Small}; this result refines the work of Pudl\'{a}k
in \cite{Pudlak-proof-lengths} and \cite{Pudlak-survey} in which \textit{%
logical depth} (length of the longest branch in the formation tree of the
formula) is used as a measure of complexity instead of the depth of
quantifier alternations complexity. Note that part (b) of the fact below is
an immediate consequence of part (a).\medskip

\noindent \textbf{Fact F.}~\textit{Suppose }$T$\textit{\ is a sequential
theory }$T$ \textit{formulated in a finite language} $\mathcal{L}$, \textit{%
and fix }$n\in \omega .$ \textit{Fix some interpretation} $\mathcal{N}$ 
\textit{of arithmetic in} $T$ \textit{satisfying }I$\Delta _{0}$.\footnote{%
It is well-known $\mathsf{Q}$ has a definable cut that satisfies I$\Delta
_{0}$ (see Theorem 5.7 of \cite{Hajek-Pudlak}), so every sequential theory
supports such an interpretation $\mathcal{N}$.}\medskip

\noindent (a) \textit{There is a} $T$-\textit{provable definable cut} $I_{n}$
of $\mathcal{N}$ \textit{and a formula} $\mathsf{Sat}_{n}(x,y)$ \textit{such
that, provably in} $T$, $\mathsf{Sat}_{n}$ \textit{satisfies the Tarskian
compositional clauses} \textit{for} $\Sigma _{n}^{\ast }$-\textit{formulae in%
} $I_{n}$ (\textit{and for all variable assignments}).\medskip

\noindent (b) \textit{There is a formula} $\mathsf{True}_{n}(x)$ \textit{%
such that, provably in T, }$\mathsf{True}_{n}(x)$ \textit{is extensional}%
\footnote{%
Without this extensionality stipulation, the numeral does not necessarily
work as a term.},\textit{\ i.e., it respects the equivalence relation
representing equality} \textit{in the interpretation} $\mathcal{N}$; \textit{%
and} \textit{for all models} $\mathcal{M}\models T$,\textit{\ and for all }$%
\Sigma _{n}^{\ast }$\textit{-sentences} $\psi ,$ \textit{we have}:

\begin{center}
$\mathcal{M}\models \left( \psi \leftrightarrow \mathsf{True}_{n}(\ulcorner
\psi \urcorner )\right) .$
\end{center}

\noindent \textbf{Proof of Theorem A.}~Suppose not, and let $T$ be
consistent completion of sequential theory formulated in a finite language $%
\mathcal{L}$. Then by the definition of sequentiality $T$ is also
sequential. Suppose to the contrary that for some $n\in \omega ,$ $T$ is
axiomatized by a set of $\Sigma _{n}^{\ast }$ sentences, i.e., suppose (1)
below:\medskip

\noindent (1) For some $n\in \omega ,$ there is a set $A$ of $\Sigma
_{n}^{\ast }$ sentences such that for all $\mathcal{L}$-sentences $\psi $, $%
\psi \in T$ iff $A\vdash \psi .$\medskip

\noindent Our proof by contradiction of Theorem A will be complete once we
verify Claim $\heartsuit $ below since it contradicts Tarski's venerable
Undefinability of Truth theorem. \medskip

\noindent CLAIM $\heartsuit $. There is a unary $\mathcal{L}$-formula $%
\varphi (x)$ such that for all $\mathcal{L}$-sentences $\psi $, $T\vdash
\psi \leftrightarrow \varphi (\ulcorner \psi \urcorner ).$\medskip

\noindent Since $T$ is sequential, we can find an $\mathcal{L}$-formula,
denoted $\mathsf{Prf}_{\mathrm{True}_{n}}(\pi ,x)$, that expresses
\textquotedblleft $\pi $ is (the code for) a first order proof of $x$ from
assumptions in $\mathsf{True}_{n}$ (i.e., from $\{x:\mathsf{True}_{n}(x)\})$%
. In particular, for each \textit{standard} $\mathcal{L}$-sentence $\psi $
and \textit{standard} $\pi $, and each model $\mathcal{M}$ of $T$, we
have:\medskip

\noindent (2) $\mathcal{M}\models \mathsf{Prf}_{\text{\textrm{True}}%
_{n}}(\pi ,\ulcorner \psi \urcorner )$ iff $\pi $ is (a code for) a proof of 
$\psi $ from $\mathsf{True}_{n}^{\mathcal{M}}:=\{\varphi :M\models \mathsf{%
True}_{n}(\ulcorner \varphi \urcorner )\}.$\medskip

\noindent Our proposed candidate of $\varphi (x)$ for establishing Claim $%
\heartsuit $ is the following formula $\rho (x)$; our choice of the letter $%
\rho $ indicates the fact that the formula expresses \textit{%
Rosser-provability }(from the true $\Sigma _{n}^{\ast }$ sentences).

\begin{center}
$\rho (x):=\exists y\left[ \mathsf{Prf}_{\mathsf{True}_{n}}(y,x)\wedge
\forall z<y\ \lnot \mathsf{Prf}_{\mathsf{True}_{n}}(z,\lnot x)\right] .$
\end{center}

\noindent Thus our goal is to show that for all $\mathcal{L}$-sentences $%
\psi $, $T\vdash \psi \leftrightarrow \rho (\ulcorner \psi \urcorner ).$ By
the completeness theorem for first order logic it suffices to show that for
each model\footnote{%
The models are just a heuristic here. In fact the whole argument can be
formulated in the complete theory $T$.} $\mathcal{M}$ of $T$, $\mathcal{M}%
\models \psi \leftrightarrow \rho (\ulcorner \psi \urcorner ).$ For the rest
of the proof, let $\mathcal{M}\models T$. We will first show:\medskip

\noindent (3)\ For all $\mathcal{L}$-sentences $\psi $, $\mathcal{M}\models
\psi \rightarrow \rho (\ulcorner \psi \urcorner ).$\medskip

\noindent To show (3), assume $\psi $ holds in $\mathcal{M}$. Let $A$\ be as
in (1), and note that $A\subseteq \mathsf{True}_{n}^{\mathcal{M}}.$ By the
assumptions about $T$, there are finitely many sentences $\alpha
_{1},...,\alpha _{n}$ in $A$ such that $\left\{ \alpha _{1},...,\alpha
_{n}\right\} \vdash \psi .$ Let $\pi _{0}\in \omega $ be (the code of) a
proof of $\psi $ from $\left\{ \alpha _{1},...,\alpha _{n}\right\} .$ Thanks
to (2) we have:\medskip

\noindent (4) $\mathcal{M}\models \mathsf{Prf}_{\text{\textsf{True}}%
_{n}}(\pi _{0},\ulcorner \psi \urcorner )$.\medskip

\noindent The assumption of consistency of $T$ coupled with (2)
yields:\medskip

\noindent (5) $\mathcal{M}\models \forall z<\pi _{0}\ \lnot \mathsf{Prf}_{%
\mathsf{True}_{n}}(z,\ulcorner \psi \urcorner ).$\medskip

\noindent This concludes the proof of (3). \medskip

\noindent To complete the proof of CLAIM $\heartsuit $, we need to show that 
$\mathcal{M}\models \lnot \psi \rightarrow \lnot \rho (\ulcorner \psi
\urcorner )$ for all $\mathcal{L}$-sentences $\psi .$ For this purpose
assume $\mathcal{M}\models \lnot \psi .$ \medskip

\noindent By putting (1) and the assumption that $\mathcal{M}\models \lnot
\psi $, we conclude that there is a \textit{standard} proof $\pi _{0}$ of $%
\lnot \psi $ from $\mathsf{True}_{n}^{\mathcal{M}},$ which by (2)
implies:\medskip

\noindent (6) For some $\pi _{0}\in \omega ,$ $\mathcal{M}\models \mathsf{Prf%
}_{\mathsf{True}_{n}}(\pi _{0},\ulcorner \lnot \psi \urcorner )$.\medskip

\noindent To see that $\mathcal{M}\models \lnot \rho (\ulcorner \psi
\urcorner )$ suppose to the contrary that $\mathcal{M}\models \rho
(\ulcorner \psi \urcorner ).$ By the choice of $\rho $, this means:\medskip

\noindent (7) For some $m_{0}\in M$, $\mathcal{M}\models \mathsf{Prf}_{%
\mathsf{True}_{n}}(m,\ulcorner \psi \urcorner )\wedge \forall z<m_{0}\ \lnot 
\mathsf{Prf}_{\mathsf{True}_{n}}(z,\ulcorner \lnot \psi \urcorner ).$
\medskip

\noindent The key observation is that putting (2) with the assumption $%
\mathcal{M}\models \lnot \psi $ allows us to conclude that the $m_{0}$ in
(7) must be a \textit{nonstandard element of }$\mathcal{M}$. Thus by
standardness of $\pi _{0}$ of (6) and the ordering properties of `natural
numbers' in $\mathcal{M},$ $\mathcal{M}\models \pi _{0}<m_{0}$, which
contradicts the second conjunct of (7).\hfill $\square $\medskip

\noindent The following result is the analogue of Theorem A for sufficiently
strong arithmetical theories. Recall that the $\Sigma _{n}$-hierarchy of
formulae is the usual hierarchy of arithmetical formulae in which $\Sigma
_{0}$-formulae are defined as formulae in which all quantifiers are bounded
(recall that $\Sigma _{0}=\Pi _{0}=\Delta _{0}$ here); $\mathrm{I}\Delta
_{0} $ is the fragment of \textsf{PA}\ in which the induction scheme is
limited to $\Delta _{0}$-formulae, and $\mathsf{Exp}$ is the arithmetical
sentence that asserts that the exponential function $2^{x}$ is total (it is
well-known that the deductive closure of $\mathrm{I}\Delta _{0}$ goes well
beyond that of $\mathsf{PA}^{-}$; moreover $\mathsf{Exp}$ is not provable in 
$\mathrm{I}\Delta _{0})$.\medskip

\noindent \textbf{Theorem B.}~\textit{For each} $n\in \omega $ \textit{every
consistent extension of} \text{I}$\Delta _{0}+\mathsf{Exp}$ (\textit{in the
same language}) \textit{that is axiomatized by a set of} $\Sigma _{n}$%
\textit{-sentences is incomplete}.\footnote{%
The set-theoretical analogue of Theorem B is Theorem C below concerning the
well-known Levy hierarchy of formulae of set theory. Theorem C can be proved
with the same strategy as in the proof of Theorem A (and the first proof of
Theorem B) thanks to the availability of the relevant definable partial
satisfaction classes in $\mathsf{KP}$. Here $\mathsf{KP}$ is Kripke-Platek
set theory with the scheme of foundation limited to $\Pi _{1}^{\mathrm{Levy}%
} $-formulae (equivalently: the scheme of $\in $-induction for $\Sigma _{1}^{%
\mathrm{Levy}}$-formulae). Thus in contrast to Barwise's $\mathsf{KP}$ in 
\cite{Barwise-book}, which includes the full scheme of foundation, our
version of $\mathsf{KP}$ is finitely axiomatizable. Note that the axiom of
infinity is not among the axioms of $\mathsf{KP}$. The existence of
definable partial satisfaction classes in $\mathsf{KP}$ follows from two
facts: (1) $\mathsf{KP}$ can prove that every set is contained in a
transitive set; and (2) $\mathsf{KP}$ can define the satisfaction predicate
for all of its internal set structures (the proofs of both of these facts
can be found in Barwise's monograph \cite{Barwise-book}; the proofs therein
make it clear that only $\Pi _{1}^{\mathrm{Levy}}$-Foundation is invoked).
See also\ Theorem 2.9 of \cite{Enayat-McKenzie} (the statement of which
involves $\mathsf{KP}$ + the axiom of infinity, but the axiom of infinity is
not used in the proof). It is also worth pointing out that $\mathsf{KP}$
plus the negation of axiom of infinity is bi-interpretable with the fragment
I$\Sigma _{1}$ of $\mathsf{PA}$ (we owe this observation to Fedor Pakhomov);
indeed the two theories can be shown to be definitionally equivalent.\medskip
\par
\noindent \textbf{Theorem C.} \textit{For each} $n\in \omega $ \textit{every
consistent completion of} $\mathsf{KP}$ (\textit{in the same language}) 
\textit{that is axiomatized by a set of} $\Sigma _{n}^{\mathrm{Levy}}$%
\textit{-sentences is incomplete.}}\medskip

\noindent \textbf{Proof.}~As shown by Gaifman and Dimitracopoulos \cite%
{Gaifman-Dimitracopoulos} (see Chapter V of \cite{Hajek-Pudlak} for an
exposition) for each $n\in \omega $ there is a formula $\mathsf{Sat}_{\Sigma
_{n}}$ such that, provably in \textrm{I}$\Delta _{0}+\mathsf{Exp,}$ $\mathsf{%
Sat}_{\Sigma _{n}}$ satisfies compositional clauses for all $\Sigma _{n}$%
-formulae. In particular there is a formula $\mathsf{True}_{\Sigma _{n}}(x)$
such that for all models $\mathcal{M}$ of \textrm{I}$\Delta _{0}+\mathsf{Exp}
$, and for all $\Sigma _{n}$-sentences $\psi $, $\psi \in \mathsf{True}%
_{\Sigma _{n}}^{\mathcal{M}}$ iff $\mathcal{M}\models \psi .$ We can now
repeat the proof strategy of Theorem A with the use of $\mathsf{True}%
_{\Sigma _{n}}^{\mathcal{M}}$ instead of $\mathsf{True}_{n}^{\mathcal{M}}.$%
\medskip

Alternatively, invoke the provability of the MRDP theorem on the Diophantine
representability of computably enumerable sets in \textrm{I}$\Delta _{0}+%
\mathsf{Exp}$ (also shown in \cite{Gaifman-Dimitracopoulos}). By the
MRDP-theorem each $\Sigma _{n}$-formula is equivalent to a $\Sigma
_{n}^{\ast }$-formula in \textrm{I}$\Delta _{0}+\mathsf{Exp}$, so Theorem A
applies.\hfill $\square $\medskip

With the help of Craig's trick\footnote{%
Suppose $T$ is computably enumerable. Fix an instance of a tautology $\tau $
in the language of $T,$ and recursively define $\tau ^{0}:=\tau $ and $\tau
^{n+1}:=\tau ^{n}\wedge \tau $. Then define $T^{\mathsf{Craig}}$ as the
result of replacing each $\varphi \in T$ with $\tau ^{n}\wedge \varphi $,
where $n$ is a witness for $\varphi \in T$. It can be readily checked that $%
T^{\mathsf{Craig}}$ is computable (indeed it is primitive recursive).} to
obtain a computable axiomatization $T^{\mathsf{Craig}}$ of an arbitrary
computably enumerable theory $T$, the proof strategy of Theorem A can be
straightforwardly adapted by using $\mathsf{Prf}_{T^{\mathsf{Craig}}+\mathsf{%
True}_{n}}$ instead of $\mathsf{Prf}_{\mathsf{True}_{n}}$ to establish the
following strengthening of Theorem A:\medskip

\noindent \textbf{Theorem A}$^{+}$. \textit{Let} $T$ \textit{be a computably
enumerable sequential theory formulated in a finite language }$\mathcal{L}$%
\textit{\ and suppose }$A$\textit{\ is a collection of }$\mathcal{L}$\textit{%
-sentences such that }$A\subseteq \Sigma _{n}^{\ast }$ \textit{for some} $%
n\in \omega $ \textit{and} $T\cup A$ \textit{is consistent. Then }$T\cup A$ 
\textit{is incomplete}.\footnote{%
Note that if $A=\varnothing $, then the proof strategy of Theorem A, when
applied to the setting of Theorem A$^{+}$, goes through for all computably
enumerable consistent extensions $T$ of the Tarski-Mostowski-Robinson theory 
$\mathsf{R}$, without the assumption of sequentiality of $T.$}\medskip

Similarly, we can obtain the following strengthening of Theorem B:\medskip

\noindent \textbf{Theorem B}$^{+}$. \textit{Let} $T$ \textit{be a computably
enumerable extension of} \textrm{I}$\Delta _{0}+\mathsf{Exp}$ (\textit{in
the same language}) \textit{and suppose }$A$\textit{\ is a collection of
arithmetical sentences such that }$A\subseteq \Sigma _{n}$ \textit{for some} 
$n\in \omega $ \textit{and} $T\cup A$ \textit{is consistent. Then }$T\cup A$ 
\textit{is incomplete}.\footnote{%
Theorem C of footnote 9 also readily lends itself to an analogous
strengthening.}\medskip

\noindent \textbf{Remark 4.}~The assumption of finiteness of language cannot
be lifted in Theorem A. For example, consider the theory $U=$ $\mathsf{CT}%
_{\omega }^{-}[\mathrm{I}\Sigma _{1}]$ of $\omega $-iterated compositional
truth over $\mathrm{I}\Sigma _{1}$ (without any induction for formulae using
nonarithmetical symbols, hence the minus superscript) formulated in an
extension of the language $\mathcal{L}_{\mathsf{A}}$ of arithmetic with
infinitely many predicates $\{\mathsf{T}_{n+1}:n\in \omega \}$, and
Tarski-style compositional axioms that stipulate that $\mathsf{T}_{n+1}$ is
compositional for all $\mathcal{L}_{n}$-formulae, with $\mathcal{L}_{0}=%
\mathcal{L}_{\mathsf{A}}$ and $\mathcal{L}_{n+1}=\mathcal{L}_{n}\cup \{%
\mathsf{T}_{n+1}\}.$ Then since bi-conditionals of form $\varphi
\longleftrightarrow \mathsf{T}_{n+1}(\ulcorner \varphi \urcorner )$ are
provable in $U$ for every $\mathcal{L}_{n}$-sentence (thanks to the
available composition axioms) ANY complete extension $V$ of $U$ is
axiomatized by $U$ (which is of bounded complexity) together with atomic
sentences of form $\mathsf{T}_{n+1}(\ulcorner \varphi \urcorner )$ where $%
\varphi \in V$ and $\varphi $ is an $\mathcal{L}_{n}$-sentence, thus $U$
axiomatizable by a set of axioms of bounded quantifier complexity. Note that 
$U$ is a sequential theory since it is an extension of the sequential theory 
$\mathrm{I}\Sigma _{1}$. By strengthening $U$ with a single axiom, namely
the so-called axiom $\mathsf{Int}$ of internal induction that uses $\mathsf{T%
}_{1}$ to state that all arithmetical instances of induction hold, we obtain 
\textit{\ a theory in an infinite language, and axiomatized by a set of
sentences of bounded quantifier complexity, whose deductive closure extends} 
$\mathsf{PA}$ \textit{and every completion of which is axiomatizable by a
set of sentences of bounded complexity}. Let us also note that it is
well-known that the technique of `$\mathcal{M}$-logic' of
Krajewski-Kotlarski-Lachlan (see Kaye's text \cite{Kaye's text} for an
exposition), or the more recent robust model-theoretic technique introduced
in \cite{Enayat-Visser} allow one to show that $\mathsf{CT}_{\omega }^{-}[%
\mathrm{I}\Sigma _{1}]+\mathsf{Int}$ is conservative over $\mathsf{PA.}$%
\medskip

Alternatively, by starting with any theory $T$ formulated in a language $%
\mathcal{L}$, we can apply a process known in model theory as \textit{%
Morleyization}\footnote{%
According to Hodges, atomization was introduced by Skolem in the 1920s, and
has \textquotedblleft nothing to do with Morley\textquotedblright . It is
classically known that the atomization of a theory can be axiomatized by
sentences of the form $\forall x_{1}...\forall x_{m}\exists y_{1}...\exists
y_{k}\ \delta $, where $\delta $ is quantifier-free. See pp.\ 62-64
(especially Theorem 2.6.6) of Hodges' majestic text \cite{Hodges-book}.} or 
\textit{atomization} to obtain an extension $T^{+}$ of $T$, formulated in an
extension $\mathcal{L}^{+}$ of $\mathcal{L}$, such that $T^{+}$ is
axiomatized by adding a collection of sentences of bounded quantifier depth
to $T$, and $T^{+}$ has elimination of quantifiers in the sense that for
each $\mathcal{L}^{+}$-formula $\varphi (x_{1},...,x_{n})$, there is an $n$%
-ary predicate $P_{\varphi }\in \mathcal{L}^{+}$ such that the equivalence $%
\varphi (x_{1},...,x_{n})\leftrightarrow P_{\varphi }(x_{1},...,x_{n})$ is
provable in $T^{+}$. The advantage of this second construction is that it
does not require the resources to build conservative truth predicates. The
atomization of a theory is well-known to be model-theoretically
conservative, whereas a truth-predicate of the type $\mathsf{T}_{1}$ in the
previous example is already not model-theoretically conservative since it
imposes recursive saturation on models of arithmetic supporting it (by a
remarkable theorem of Lachlan; see Kaye's text \cite{Kaye's text} for an
exposition).\medskip

\noindent \textbf{Remark 5.}~In Theorem B, the theory \textrm{I}$\Delta _{0}+%
\mathsf{Exp}$ cannot be weakened to $\mathsf{PA}^{-}$, i.e., for some $n\in
\omega $ there is a consistent completion of $\mathsf{PA}^{-}$ (in the same
language) that is\textit{\ }axiomatized by a set of $\Sigma _{n}$\textit{-}%
sentences.\textit{\ }The proof of this and related results will appear in
our upcoming paper with Mateusz \L e\l yk. It will be hard to prove an
analogous result about $\mathrm{I}\Delta _{0}$ since the proof strategy of
Theorem B makes it clear that the analogous result for $\mathrm{I}\Delta
_{0} $ implies that $\mathrm{I}\Delta _{0}$ does not prove that there is a $%
\Sigma _{1}$-satisfaction predicate, which is known to be a tall order, as
indicated in \cite{Zossia+Leszek+Jeff}.\medskip

\noindent \textbf{Remark 6.}~The proof of Theorem A/A$^{+}$ can be used to
establish:\medskip

\noindent $(\triangledown )$ \textit{There is no consistent complete theory }%
$T$\textit{\ interpreting }$\mathsf{R}$ \textit{that is of the form} $%
T=T_{0}+A$, \textit{where} $T_{0}$ \textit{is computably enumerable, and} $A$%
\textit{\ is a subset of some decidable set} $\Gamma $ \textit{of sentences
for which }$T$\textit{\ has a truth predicate. }\medskip

\noindent Theorem A/A$^{+}$ then follows by putting $(\triangledown )$, with
the choice of $\Gamma :=\Sigma _{n}^{\ast }$, together with Fact F.
Similarly, Theorems B/B$^{+}$ and Theorem C can be obtained from $%
(\triangledown )$ by invoking the corresponding well-known analogues of Fact
F. This more general approach is further explored in \cite{Visser-parallel
paper}, but it should be noted that upon the completion of our paper we were
kindly informed by Emil Je\v{r}\'{a}bek that $(\triangledown )$ was
formulated and established in his MathOveflow answer \cite{Emil-MO}.

\medskip

We conclude with a question that is motivated by Theorem A and the fact that 
$\mathsf{Q}$ is not a sequential theory \cite{Visser on Q}.\medskip

\noindent \textbf{Question }$\#$\textbf{.}~\textit{Is it possible for a
consistent completion of} $\mathsf{Q}$ \textit{to be axiomatized by a
collection of sentences of bounded quantifier-depth}?

\noindent \textsc{Ali Enayat, Department of Philosophy, Linguistics, and
Theory of Science, University of Gothenburg}, \textsc{\ Sweden.}

\noindent \texttt{ali.enayat@gu.se}\medskip

\noindent \textsc{Albert Visser, Department of Philosophy, Faculty of
Humanities, Utrecht University, The Netherlands.}

\noindent \texttt{albert.visser@uu.nl}

\end{document}